\input amstex \documentstyle{amsppt}\magnification=1200 \nologo
\NoBlackBoxes  \hsize=16 truecm \vsize = 9 truein \voffset = 7.5
truemm \hoffset = 0.1 truemm  \topmatter

\title
Canonical map of low codimensional subvarieties.
\endtitle

\author Valentina Beorchia*, Ciro Ciliberto** and Vincenzo Di Gennaro*** \endauthor

\leftheadtext{Valentina Beorchia, Ciro Ciliberto and Vincenzo Di
Gennaro} \rightheadtext{Canonical map of low codimensional
subvarieties}

\medskip
\abstract \nofrills

ABSTRACT. Fix integers $a\geq 1$, $b$ and $c$. We prove that for
certain projective varieties $V\subset{\bold P}^r$ (e.g. certain
possibly singular complete intersections), there are only finitely
many components of the Hilbert scheme parametrizing irreducible,
smooth, projective, low codimensional subvarieties $X$ of $V$ such
that
$$
h^0(X,\Cal O_X(aK_X-bH_X)) \leq \lambda
d^{\epsilon_1}+c\left(\sum_{1\leq h <
\epsilon_2}p_g(X^{(h)})\right),
$$
where $d$, $K_X$ and $H_X$ denote the degree, the canonical
divisor and the general hyperplane section of $X$,  $p_g(X^{(h)})$
denotes the geometric genus of the general linear section of $X$
of dimension $h$, and where $\lambda$, $\epsilon_1$ and
$\epsilon_2$ are suitable positive real numbers depending only on
the dimension of $X$, on $a$ and on the ambient variety $V$. In
particular, except for finitely many families of varieties, the
canonical map of any irreducible, smooth, projective, low
codimensional subvariety $X$ of $V$, is birational.



\endabstract
\endtopmatter

{\footnote""{*Dipartimento di Scienze Matematiche,
Universit\`a\ di Trieste, via Valerio 12, 34100 Trieste, Italy.
E-mail: beorchia{\@}units.it}

{\footnote""{**Dipartimento di Matematica, Universit\`a\ di Roma
\lq \lq Tor Vergata", 00133 Roma, Italy. E-mail:
cilibert{\@}axp.mat.uniroma2.it}

{\footnote""{ ***Dipartimento di Matematica, Universit\`a\ di Roma
\lq \lq Tor Vergata", 00133 Roma, Italy. E-mail:
digennar{\@}axp.mat.uniroma2.it}}

\bigskip

\flushpar{\bf{0. Introduction.}} \bigskip

A famous theorem of Ellingsrud and Peskine [EP] states that there
are only finitely many components of the Hilbert scheme
parametrizing smooth surfaces in ${\bold P}^4$ not of general
type. This paper has been followed by others in which suitable
extensions have been presented (see [AS], [DC] , [BOSS], [FO],
[S]). More recently two of us gave further wide extensions of
these results ([CD1], [CD2], [CD3]).

Going back to the original theorem of Ellingsrud and Peskine,
Ellia and Folegatti [EF] remarked that the technique of proof
makes it possible to show a more general result. Namely they are
able to prove boundedness for families of smooth surfaces in
${\bold P}^4$ with geometric genus bounded above by the sectional
genus. This in turn implies boundedness for families of smooth
surfaces in  ${\bold P}^4$ with non birational canonical map.\par

The present paper is devoted to give a wide extension of
Ellia-Folegatti's result (see Theorem 0.1), which we now will
state.

Let $V$ be an irreducible, possibly singular, projective variety
over $\bold C$. Let ${\Cal S}$ be a set of projective subvarieties
of $V$. We will say that ${\Cal S}$ is {\it bounded} if there is a
closed immersion $V\subset \bold P^r$ such that
$$
\sup\{deg(X): X\in {\Cal S}\}<+\infty.
$$
This means that the varieties in ${\Cal S}$ belong to finitely
many components of the Hilbert scheme. In particular, this
definition does not depend on the closed immersion.

In this paper we will prove the following:
\medskip

\proclaim {Theorem 0.1} Let $V\subset {\bold P^r}$  be an
irreducible, projective variety of dimension $m$. Let $1\leq n<m$
be an integer, and put $k=m-n$. Fix integers $a,\,b,\,c \in {\bold
Z}$, with $a\geq 1$, and put $\epsilon(a)=min \lbrace
\frac{n}{k}+1,\, \frac{n}{k}+a-1\rbrace$.  Assume that at least
one of the following properties holds.

\medskip
(A) $m=n+2$, $2\leq n\leq 4$, $V$ is smooth, $NS(V)\simeq{\bold
Z}$ and any algebraic class in $H^{4n-8}(V,{\bold C})$ is a
multiple of $H^{2n-4}_V$, where $H_V$ is a hyperplane section of
$V$;

\medskip
(B) $m=n+2$, $n\geq 4$ and, only when $a=1$, $n\geq
\frac{r+1}{2}$; for $i=-1,0$ any algebraic class in
$H_{2n+2i}(V,{\bold C})$ is a multiple of $H^{2-i}_V$, the general
linear section $V^{(4)}$ of dimension $4$ of $V$ is smooth and
$NS(V^{(4)})\simeq{\bold Z}$;

\medskip
(C) $m=n+2$, $n\geq 5$, $V$ is smooth, for $1\leq i\leq 3$ any
algebraic class in $H^{2i}(V,{\bold C})$ is a multiple of
$H^{i}_V$, and there exist rational numbers $v_1,\dots,v_n$ such
that $c_i(T_V)=v_iH_V^i$ in $H^{2i}(V,{\bold C})$ for any $1\leq
i\leq n$;

\medskip
(D) $ n\geq \frac{m+2}{2}$ and, only when $a=1$, $n\geq
\frac{r+1}{2}$; for $i=-1,0$ any algebraic class in
$H_{2n+2i}(V,{\bold C})$ is a multiple of $H^{k-i}_V$; moreover,
for some $2\leq h \leq n$ with $h\geq k$, the general linear
section $V^{(h+k)}$ of dimension $h+k$  of $V$ is smooth, and
either  any algebraic class in $H^{2i}(V^{(h+k)},{\bold C})$ is a
multiple of $H^i_{V^{(h+k)}}$ for $i\in \lbrace 1,\, k\rbrace$, or
$k$ is even and $H^{2i}(V^{(h+k)},{\bold C})\simeq \bold C$ for
any $i=1,\dots,k-1$;

\medskip
(E) $n > \frac{3m-2}{4}$, $V$ is smooth, $a\geq 2$,  for $1\leq
i\leq 2k-1$ any algebraic class in $H^{2i}(V,{\bold C})$ is a
multiple of $H^{i}_V$, and there exist rational numbers
$v_1,\dots,v_n$ such that $c_i(T_V)=v_iH_V^i$ in $H^{2i}(V,{\bold
C})$ for any $1\leq i\leq n$.

\medskip
For any $n$-dimensional subvariety $X$ of $V$ denote by $X^{(h)}$
($K_X$ resp.) the general linear section of dimension $h$ (the
canonical divisor resp.) of $X$. Denote by $p_g(X^{(h)})$ the
geometric genus of $X^{(h)}$. Put $H_X= X^{(n-1)}$ and $d=deg(X)$.

Then there exists a strictly positive real number $\lambda
>0$, depending only on $n$, $a$ and the ambient variety $V$, such
that the set of irreducible, smooth, projective subvarieties $X$
of $V$ of dimension $n$ satisfying the following inequality under
the hypotheses (A) with $n\neq 3$, or (C) or (E)
$$
h^0(X,\Cal O_X(aK_X-bH_X)) \leq \lambda d^{{\frac{n}{k}}+1}+
c\left( \sum_{h=1}^{n-1} p_g(X^{(h)})\right), \tag 0.1.1
$$
or the following inequality under the hypotheses (A) with $n= 3$
$$
h^0(X,\Cal O_X(aK_X-bH_X)) \leq \lambda
d^2+c(p_g(X^{(1)})+p_g(X^{(2)})),
$$
or the following inequality under the hypotheses (B) or (D)
$$
h^0(X,\Cal O_X(aK_X-bH_X)) \leq \lambda d^{\epsilon(a)}+ c\left(
\sum_{1\leq h < {{\epsilon (a)}-1}} p_g(X^{(h)})\right), \tag
0.1.2
$$
is bounded.

Moreover, when $V$ is smooth, under the hypotheses (A) with $n\neq
3$, or  (B) with $a\geq 2$, or (C), or  (D) with $a\geq 2$, or
(E), the previous estimates are sharp in the following sense:
there exists a real number $\mu >\lambda$ depending only on $n$,
$a$ and the ambient variety $V$, such that the set of irreducible,
smooth, projective subvarieties $X$ of $V$ of dimension $n$
satisfying the inequality
$$
h^0(X,\Cal O_X(aK_X-bH_X)) \leq \mu d^{{\frac{n}{k}}+1} \tag 0.1.3
$$
is not bounded.
\endproclaim \bigskip

As a consequence one has the:

\medskip
\proclaim {Corollary 0.2} With the same assumption of  Theorem
0.1, except for finitely many families of varieties, the canonical
map of any irreducible, smooth, projective subvariety of $V$ of
dimension $n$, is birational.
\endproclaim \medskip

Theorem 0.1 has a rather wide range of applications, also to
singular varieties. By Lefschetz Hyperplane Theorem,
Poincar\graveaccent e duality and Barth Theorem, any smooth
complete intersection fourfold $V$ on a Grassmann variety or on a
spinor variety [PR], any smooth complete intersection $V\subset
{\bold P^r}$ of dimension $5$ or $6$, any smooth fourfold in
${\bold P^6}$ and any smooth sixfold in ${\bold P^8}$ verifies the
hypothesis (A).

When $n\geq 5$ (and, only for $a=1$, when $n\geq {(r+1)}/{2}$),
possibly singular complete intersections $V\subset {\bold P^r}$ of
dimension $n+2$, with $dim (Sing(V))\leq n-6$, verify the
hypothesis (B) (see [Dm], Theorem $(2.11)$, pg.144). Moreover any
hypersurface $V\subset {\bold P^{n+3}}$ of degree $t\geq 3$
defined by the equation
$$
x_0^ax_1^{t-a}+x_1x_2^{t-1}+\dots+x_{n+1}x_{n+2}^{t-1}+x_{n+3}^{t}=0,
$$
with $n\geq 4$, $1\leq a < t-1$ and $(a,t)=1$, has at most two
singular points and satisfies the hypothesis (B) (see [Dm],
Proposition $(2.24)$, pg.148). This provides examples for the
assumption (B) in the case $n=4$ too. By Barth Theorem again, any
smooth subvariety $V\subset {\bold P^r}$ of dimension $n+2$, with
$n\geq (r+2)/2$, satisfies the hypothesis (B).

As before one sees that any possibly singular complete
intersection $V\subset {\bold P^r}$ of dimension $m$ with  $n >
(m+2)/2$ and $dim(Sing(V))\leq 2n-m-4$ (and, only for $a=1$, with
$n\geq {(r+1)}/{2}$), and any smooth subvariety $V\subset {\bold
P^r}$ of dimension $m$ with $n\geq (r+2)/2$, satisfies the
hypothesis (D). When $n=(m+2)/2$ (and, only for $a=1$, when $n\geq
{(r+1)}/{2}$), by Noether-Lefschetz Theorem [CGH] this is true
also for Noether-Lefschetz general complete intersections
$V\subset {\bold P^r}$ of dimension $m$, with $h^{a,m-a}(V)\neq 0$
for some $a<m/2$ (for instance, any general hypersurface of  even
dimension $m\geq 4$ and degree $\geq 3$ [PS]).

Smooth complete intersections $V\subset {\bold P^r}$ of dimension
$m\geq 7$ ($m<(4n+2)/3$ resp.) satisfy the hypothesis (C) ((E)
resp.).
\smallskip

As for the proof, the general strategy consists in determining a
lower bound for the geometric genus of a subvariety $X\subset V$.
Our methods are  partly based on the technical developments of
[CD1], [CD2] and [CD3], to which we will often refer.

Under the hypothesis (A), the basic tools are inequality $(1.5)$
proved in [CD1], Castelnuovo-Halphen's theory, and Miyaoka-Yau's
inequality (see Section 1).

The line of the proof under the hypotheses (B) is the following. A
Barth-Lefschetz type of argument proves that the N\'eron-Severi
group of the subvarieties $X$ of $V$ has rank 1. Then, when $a\geq
2$, Kawamata-Viehweg Vanishing Theorem allows us to obtain a lower
bound for $h^0(X,\Cal O_X(aK_X-bH_X))$ (see Proposition 2.2, which
should be compared with Koll\'ar-Luo-Matsusaka estimate [K], pg.
302, Theorem 2.15.9). In order to make this lower bound explicit,
and then deduce the boundedness for $X$, we need a general result,
i.e. Theorem 2.1, concerning boundedness for subvarieties with
bounded sectional genus, whose proof relies on inequality $(1.5)$
and Castelnuovo theory, and does not need the assumption $n\geq
(r+1)/2$. The previous argument does not work when $a=1$. In this
case, using the hypothesis $n\geq (r+1)/2$, we may apply Larsen
Theorem (as in Amerik paper [A]) and deduce that the Picard group
of $X$ is generated by the hyperplane section. Then Castelnuovo
theory  [EH] enables us to bound $h^0(X,\Cal O_X(K_X-bH_X))$ from
below, and so, using Theorem 2.1 again,  one may conclude in a
similar way as in the case $a\geq 2$ (see Section 2).

Under the hypotheses (C), the Chern classes of the normal bundle
of a smooth subvariety $X\subset V$ are multiples of the linear
sections. By means of a somewhat delicate numerical analysis based
on the Hirzebruch-Riemann-Roch Theorem and the previous Theorem
2.1, this allows us to  bound from below the arithmetic genus of
$X$ in terms of the degree $d$ of $X$. Using the Hyperplane
Lefschetz Theorem one may estimate the difference between the
geometric genus of $X$ and the arithmetic genus, finally obtaining
a lower bound for the geometric genus, from which one easily
concludes (see Section 3).

Under the hypotheses (D), Theorem 0.1 follows in a similar manner
as under the hypotheses (B), taking into account a general result,
i.e. Theorem 4.1, which states boundedness for subcanonical
subvarieties, and does not need the assumption $n\geq (r+1)/2$.
For the proof of this result, the main tools are Chern classes
computations like in [S] and [CD3] (see Section 4).

Under the hypotheses (E), Theorem 0.1 follows by combining the
methods used in the proof under the hypotheses (C) and  (D) (see
Section 5).

\smallskip

We  prove the sharpness of the estimates $(0.1.1)$ and $(0.1.2)$
in the sense of $(0.1.3)$, by considering the unbounded set
consisting of the $k$-codimensional  complete intersections of
balanced type on $V$.

For instance, for a smooth complete intersection surface  $X$ of
balanced type $(u,u)$ in ${\bold P^4}$, we have $p_g(X)\leq
{\frac{7}{12}}d^2$. On the other hand, our analysis proves that
smooth surfaces $X\subset {\bold P^4}$ with $p_g(X)\leq
{\frac{1}{12}}d^2$ are bounded (see Section 1, $(1.2)$).

Under the hypotheses (A) with $n=3$, or (B) with $a=1$, or (D)
with $a=1$, our methods do not enable us to obtain the expected
sharp estimate. In any case we give explicit estimate for the
constant $\lambda$, and our analysis may give, in principle,
explicit bounds for the degree $d$ in terms of the given data. We
decided not to dwell on this here.

It would be interesting to investigate whether Theorem 0.1 is
sharp in a stronger sense. For instance, it implies boundedness
for smooth surfaces of degree $d$ in ${\bold P^4}$ with geometric
genus $\leq \lambda d^2$, with $\lambda <1/6$ (see Section  1, in
particular $(1.2)$). A nice question is therefore whether the
family of smooth surfaces in ${\bold P^4}$ with geometric genus
$\leq d^2/6$ is bounded or not.

\bigskip
\bigskip
\bigskip

\noindent {\bf {Notation.}}\bigskip Let $Y$ be any smooth,
irreducible, projective variety over ${\bold C}$. We will denote
by $T_Y$ the tangent bundle of $Y$, and by $K_Y$ a canonical
divisor of $Y$. If ${\Cal E}$ is any sheaf on $Y$ we denote by
$\chi({\Cal E})$ its Euler-Poincar\'e characteristic and by
$c_i({\Cal E})$ its Chern classes. We denote by $p_g(Y)$  the
geometric genus of $Y$. As usual $NS(Y)$ will be the
N\'eron-Severi group of $Y$. We denote numerical equivalence by
using the symbol  $\equiv$ . If $Z\subset Y$ is a subvariety, we
denote by $N_{Z,Y}$ the normal sheaf  of $Z$ in $Y$. When
$Y\subset{\bold P^{r}}$ we denote by $H_Y$ the general hyperplane
section of $Y$. Moreover, if $dim(Y)=l$ and $0\leq j \leq l$, we
denote by $Y^{(j)}$ the intersection of $Y$ with a general linear
subspace $\bold P^{r-l+j}\subset\bold P^{r}$. In particular
$dim(Y^{(j)})=j$, and $Y^{(l-1)}=H_Y$. We say that $Y\subset{\bold
P^{r}}$ is {\it {numerically subcanonical}}  if $K_Y=eH_Y$ in
$H^{2}(Y,{\bold Q})$ for some $e\in {\bold Q}$. We say that
$Y\subset{\bold P^{r}}$ is {\it{subcanonical}} if $K_Y=eH_Y$ in
$Pic(Y)$, for some $e\in {\bold Z}$.

If $x$ is a real number, we denote by $[x]$ the integral part of
$x$. If ${\Cal S}$ is a set and $f: {\Cal S}\to {[0,+\infty[}$ is
a numerical function, we say that a function $\phi: {\Cal S}\to
{[0,+\infty[}$ is $O(f)$, and we write $\phi=O(f)$, if
$|\phi(\xi)|\leq C f(\xi)$ for all $\xi\in {\Cal S}$, where $C$ is
a constant $>0$.

\bigskip
\bigskip
\bigskip

\noindent {\bf 1. The proof of Theorem 0.1 under the hypotheses
(A).}
\bigskip

We start by proving  Theorem 0.1 in the case $n=2$, under the
hypotheses (A). We need the following result:
\medskip
\proclaim {Theorem 1.1} Let $V\subset {\bold P^r}$  be an
irreducible, projective variety of dimension $n+2\geq 4$. Assume
that the general linear section $V^{(4)}$ of dimension $4$  of $V$
is smooth and such that $NS(V^{(4)})\simeq{\bold Z}$. Fix an
integer $s$, and a real number
$$
\lambda < {\frac {1}{(n+1)!s^n}}.
$$
For any projective $n$-fold $X$ contained in $V$, put $d=deg(X)$
and let $p_g(X)$ be the geometric genus of $X$. Then the set of
irreducible, smooth, projective, codimension two subvarieties $X$
of $V$,  contained in some reduced, projective subvariety of
${\bold P^r}$ of dimension $n+1$ and degree $\leq s$, and such
that
$$
p_g(X)\leq \lambda d^{n+1}, \tag 1.1.1
$$
is bounded.
\endproclaim \medskip

\demo {Proof of Theorem 1.1} Using the same argument as in the
proof of [CD2], Theorem 4.1, pg. 487, one proves that
$$
p_g(X)={\frac {d^{n+1}}{(n+1)!s^n}}+O(d^n) \tag 1.1.2
$$
for any irreducible, smooth, projective, codimension two
subvarieties $X$ of $V$ of degree $d$,  contained in some reduced,
projective subvariety of ${\bold P^r}$ of dimension $n+1$ and
degree $s$ (see [CD2], pg. 491, line 7 from below). Our Theorem
$(1.1)$ follows comparing $(1.1.1)$ with $(1.1.2)$. $\square$
\bigskip

We are in position to prove  Theorem 0.1 in the case $n=2$. To
this purpose, let $V\subset {\bold P^r}$ be a smooth fourfold as
in  Theorem 0.1. We may assume $V$ is nondegenerate. First we
examine the case $a=1$ and $b=0$. Put $t= deg(V)$. Fix any real
number $\lambda$ such that
$$
\lambda <  {\frac{1}{6t}}. \tag 1.2
$$
Let $X\subset V$ be any smooth projective surface such that
$$
h^0(X,\Cal O_X(K_X))\leq \lambda d^2+cp_g(H_X). \tag 1.3
$$
We have
$$
\chi(\Cal O_X)\leq 1+\lambda d^2+cp_g(H_X). \tag 1.4
$$
Now we use the crucial inequality (see [CD1], pg.277, (2))
$$
d^2/t+q(2g-2)+O(d)\leq 2(K^2_X-6\chi(\Cal O_X)), \tag 1.5
$$
where $g=p_g(H_X)$ and $K_V=qH_V$ in $H^{2}(V,{\bold C})$. By
[CD1] we  may assume that $X$ is of general type. Hence, using
Miyaoka-Yau  inequality (compare with [EF])
$$
K^2_X\leq 9\chi(\Cal O_X)
$$
and $(1.4)$ and $(1.5)$, we get
$$
\bigg({\frac{1}{t}}-6\lambda\bigg)d^2 +(2q-6c)g+O(d)\leq 0. \tag
1.6
$$
If $2q-6c\geq 0$ then $d$ is bounded by $(1.2)$. So we only have
to examine the case $2q-6c< 0$. From $(1.6)$ we get
$$
g\geq d^2\frac {1-6t\lambda}{t(6c-2q)}+O(d). \tag 1.7
$$
Put
$$
\frac {1}{2\alpha}=\frac {1-6t\lambda}{t(6c-2q)},
$$
and
$$
s_0=max \lbrace t+1,\, [\alpha]+1\rbrace.
$$
Notice that, by $(1.3)$ and Castelnuovo's bound on the geometric
genus of a projective curve, we have
$$
p_g(X)\leq \lambda d^2+cp_g(H_X)=O(d^2).
$$
Then, by  Theorem 1.1 we may assume that there is no irreducible
and reduced 3-fold $T\subset {\bold P^r}$ containing $X$ and of
degree $<s_0$. Hence, by [CC], $H_X$ is nondegenerate and not
contained in any surface in ${\bold P^{r-1}}$ of degree $<s_0$. By
[CCD] we deduce
$$
g\leq d^2/2s_0+O(d). \tag 1.8
$$
From $(1.7)$ and $(1.8)$ we obtain
$$
\frac {d^2}{2}\left(\frac {1}{\alpha}-\frac
{1}{s_0}\right)+O(d)\leq 0
$$
which implies that  $d$ is bounded because $s_0>\alpha$. This
concludes the proof of the boundedness in Theorem 0.1, under the
hypotheses (A), with $n=2$, $a=1$ and $b=0$.

Now we turn to the case $n=2$ under the hypotheses (A), with $a=1$
and any fixed $b$. As before, fix any real number $\lambda$ as in
$(1.2)$. Let $X\subset V$ be a smooth surface such that
$h^0(X,\Cal O_X(K_X-bH_X))\leq \lambda d^2+cp_g(H_X)$. From the
Poincar\graveaccent e residue sequence (compare with [E], proof of
Corollary $(2.2)$ (a), and with [BOSS], pg. 329)
$$
0\to \omega_X(-1)\to \omega_X\to \omega_{H_X}(-1)\to 0
$$
we deduce
$$
p_g(X)\leq \lambda d^2+(b+c)p_g(H_X).
$$
And so the boundedness of $d$ follows from the previous analysis
of the case $a=1$, $b=0$.

Next, we consider the case $a\geq 1$.  Fix any real number
$\lambda$ as in $(1.2)$, and let $X\subset V$ be a smooth surface
such that $h^0(X,\Cal O_X(aK_X-bH_X))\leq \lambda d^2+cp_g(H_X)$.
Again by the previous analysis, we may assume $p_g(X)>0$. Hence we
have
$$
h^0(X,\Cal O_X(K_X-bH_X))\leq h^0(X,\Cal O_X(aK_X-bH_X)).
$$
Therefore the boundedness of $d$ follows from the analysis of the
case $a=1$.
\smallskip

Finally, to prove the sharpness of  estimate $(0.1.1)$ under the
hypotheses (A) with $n =2$, consider the set $\Cal S$ of  smooth
surfaces complete intersection on $V$ of type $(u,u)$. Clearly,
$\Cal S$ is not bounded. In order to estimate  $h^0(X,\Cal
O_X(aK_X-bH_X))$ for $X\in \Cal S$, first notice that $d=tu^2$
($t$ = degree of $V$), and by the adjunction formula we have
$K_X=(q+2u)H_X$ in $H^{2}(X,{\bold C})$, where $K_V=qH_V$ in
$H^{2}(V,{\bold C})$. From [K], pg. 301, $(2.15.8.6)$, we have
$$
h^0(X,\Cal O_X(aK_X-bH_X))\leq
{\frac{\left[(a(q+2u)-b)H_X^2\right]^2}{H_X^2}}+2={\frac{4a^2d^2}{t}}+
O(d^{\frac{3}{2}}).
$$
This proves $(0.1.3)$ and completes the proof of Theorem 0.1 under
the hypotheses (A) with $n =2$.
\bigskip \bigskip

Next we are going to prove  Theorem 0.1 in the case $n=3$, under
the hypotheses (A). First we examine the case $a=1$ and $b=0$. To
this purpose, taking into account Castelnuovo's bound for the
genus of a projective curve, it suffices to prove the boundedness
of the set of smooth $3$-folds $X\subset V$ such that
$$
p_g(X)\leq \lambda d^2 +cp_g(H_X), \tag 1.9
$$
where $\lambda$ is {\it any} fixed real number, and $c$ is any
fixed integer $\geq 1$. We follow the proof of Theorem 0.1 in the
case $n=3$, in [CD2]. First notice that using  $(2.2)$ of [CD2]
and Hyperplane Lefschetz Theorem one has
$$
\multline \chi (\Cal O_X )\geq 1-h^1(X,\Cal O_X)-p_g(X)\\ \geq
c(1-h^1(X,\Cal
O_X)-p_g(H_X))+(cp_g(H_X)-p_g(X)+1-c)\\
\geq c(-\chi (\Cal
O_{H_X})-2(g-1))+(cp_g(H_X)-p_g(X)+1-c),\\
\endmultline
$$
where $g=p_g(X^{(1)})$ denotes the linear sectional genus of $X$.
As in the proof of (i) of Lemma 2.1 in [CD2], we deduce
$$
\multline (24\gamma-12q+24c)\chi (\Cal O_{H_X} ) \geq
d^2(2\gamma-q-2)/t+(g-1)(2d/t+O(1))+O(d)\\
+24(cp_g(H_X)-p_g(X)),\\
\endmultline
$$
where $\gamma$ is any integer such that the twisted tangent bundle
$T_V(\gamma)$ is globally generated and $K_V=qH_V$ in
$H^2(V,{\bold C})$. Using (g) of Lemma 2.1 in [CD2], and taking
into account that we may assume $2\gamma-q+2c>0$, we obtain
$$
\multline
 0\leq
(g-1)[\mu(g-1)/d-2d/t+O(1)]-2d^2(2\gamma-q+c-1)/t+O(d)\\
+24(p_g(X)-cp_g(H_X)),\\
\endmultline
$$
where $\mu =8(2\gamma-q+2c)$. From $(1.9)$ we get
$$
0\leq (g-1)[\mu(g-1)/d-2d/t+O(1)]-2d^2[-12\lambda+
(2\gamma-q+c-1)/t]+O(d). \tag 1.10
$$
By Theorem 1.1 in [CD2]  we may assume $g>1$. Moreover, as in the
case $n=2$, using Theorem 1.1 and $(1.9)$, we may also assume
$$
g\leq d^2/2s_0+O(d),
$$
where
$$
s_0=max \lbrace t+1,\, [\mu t/4]+1\rbrace.
$$
Therefore
$$
(g-1)[\mu(g-1)/d-2d/t+O(1)]\leq 0
$$
for $d>O(1)$, and by $(1.10)$ we get
$$
2d^2[-12\lambda+ (2\gamma-q+c-1)/t]+O(d)\leq 0,
$$
which proves the boundedness of $d$, because we may choose
$\gamma$ such that $-12\lambda+ (2\gamma-q+c-1)/t >0$. This
concludes the proof of  Theorem 0.1 under the hypotheses (A), in
the case $n=3$, when $a=1$ and $b=0$.

As in the case $n=2$, one reduces the proof of the general case
$a\geq 1$ and  $b\, ,c\in {\bold Z}$, to the case $a=1$ and $b=0$.
This concludes the proof of  Theorem 0.1 under the hypotheses (A),
in the case $n=3$. \bigskip \bigskip

Now we are going to prove the Theorem in the case $n=4$, under the
hypotheses (A). First we examine the case $a=1$ and $b=0$. Fix any
positive real number $\lambda$ such that
$$
\lambda < \frac {1}{1440t^2}. \tag 1.11
$$
As before, taking into account Castelnuovo's bound for the genus
of a projective curve, it suffices to prove the boundedness of the
set of smooth projective $4$-folds $X\subset V$ such that
$$
p_g(X)\leq \lambda d^3 +c(p_g(X^{(2)})+p_g(X^{(3)})), \tag 1.12
$$
where  $c$ is any fixed integer $\geq 1$. We follow the proof of
Theorem 0.1 in the case $n=4$, in [CD2]. Arguing as above, we
obtain
$$
\multline \chi (\Cal O_X )\leq 1+p_g(X^{(2)})+p_g(X)\\
\leq
(c+1)(1+p_g(X^{(2)})+p_g(X^{(3)}))+(p_g(X)-cp_g(X^{(2)})-cp_g(X^{(3)}))\\
\leq (c+1)(2\chi (\Cal O_{X^{(2)}})+2(g-1)+2-\chi (\Cal
O_{X^{(3)}}))+(p_g(X)-cp_g(X^{(2)})-cp_g(X^{(3)})),\\
\endmultline
$$
where $g=p_g(X^{(1)})$ denotes the linear sectional genus of $X$.
As in the proof of (m) of Lemma 3.1 in [CD2], we deduce
$$
\multline 60(2\gamma-q+2(c+1))\chi (\Cal O_{X^{(3)}})\\ \leq
(-9d/t+O(1))\chi (\Cal
O_{X^{(2)}})+(g-1)[d(9q/2-9-10\gamma)/t+O(1)]\\-d^3/(12t^2)+O(d^2)
+120(p_g(X)-cp_g(X^{(2)})-cp_g(X^{(3)})),\\
\endmultline \tag 1.13
$$
where $\gamma$ is any integer such that $T_V(\gamma)$ is globally
generated and $K_V=qH_V$ in $H^2(V,{\bold C})$. Now put
$$
\alpha=2g-2-3d,\quad\quad \nu =2\gamma-q, \quad \quad \beta =\frac
{5(2\gamma-q+2(b+1))}{2(1+\nu d/\alpha)}.
$$
Notice that we may assume $\gamma >O(1)$. By Theorem 1.1 in [CD2]
we also may assume
$$
g\geq O(d)\quad {\text {and}} \quad d\geq O(1). \tag 1.14
$$
Hence we may assume $0<\nu d/\alpha<1/2$, from which we get
$$
0< 5(2\gamma -q+2(b+1))/3 <\beta <  5(2\gamma -q+2(b+1))/2. \tag
1.15
$$
Combining inequality (i) of Lemma 3.1 in [CD2], i.e.
$$
\multline 24(1+\nu d/\alpha)\chi (\Cal O_{X^{(3)}})\geq
(g-1)(2d/t+O(1))\\
+\chi (\Cal O_{X^{(2)}})(-6d^2/(\alpha t)+O(1))-36(\chi (\Cal
O_{X^{(2)}}))^2/\alpha -d^4/(4\alpha t^2)+O(d^2)\\
\endmultline
$$
(notice that the term $36(\chi (\Cal O_{X^{(2)}}))^2/\alpha$
corrects a misprint in [CD2]), with $(1.12)$ and $(1.13)$, we get
$$
A(\chi (\Cal O_{X^{(2)}}))^2+B\chi (\Cal O_{X^{(2)}})+C\geq 0,
$$
where $A=36$, $B=6d^2/t-9d\alpha/(\beta t)+\alpha O(1)$, and
$$
\multline
 C=\alpha(g-1)[d(9q/2-9-10\gamma)/(\beta t)
-2d/t+O(1)]\\
+d^4/(4t^2)+\alpha d^3(120\lambda-1/(12t^2))/\beta +\alpha
O(d^2).\\
\endmultline
$$
Using $(1.11)$, $(1.14)$ and $(1.15)$ one sees that $B<0$ and
$C<0$. At this point, taking into account Theorem 1.1, in order to
prove the boundedness of $d$, one may proceed exactly as in the
proof of Theorem 0.1 in [CD2], in the case $n=4$, under the
hypotheses (A) (see [CD2], pg. 486, line 18 from above). This
concludes the proof of  Theorem 0.1 in the case $n=4$, under the
hypotheses (A), when $a=1$, $b=0$ and $c$ is any integer. Next, as
before, one reduces the general case $a\geq 1$ and  $b\, ,c\in
{\bold Z}$, to the case $a=1$ and $b=0$.

One proves the sharpness in the sense of $(0.1.3)$ in a similar
way as in the case $n=2$.

This concludes the proof of  Theorem 0.1 in the hypotheses (A).
\bigskip
\noindent {\bf Remark 1.16.} From the previous proof of Theorem
0.1 under the hypotheses (A), we see that in the case $n=3$, for
any fixed integers $a\geq 1$, $b$, $c$, and for {\it{any}} fixed
real number $\lambda$, the set of irreducible, smooth, projective
subvarieties $X$ of $V$ of dimension $3$ such that $ h^0(X,\Cal
O_X(aK_X-bH_X)) \leq \lambda d^2+c(p_g(X^{(1)})+p_g(X^{(2)}))$, is
bounded. $\square$

\bigskip
\bigskip
\bigskip

\noindent {\bf 2. The proof of  Theorem 0.1 under the hypotheses
(B).}
\bigskip

We begin with  the following preliminary result, concerning
boundedness for subvarieties with bounded sectional genus (compare
with Theorem 1.1 in [CD2]).
\medskip

\proclaim {Theorem 2.1} Let $V\subset {\bold P^r}$  be an
irreducible, projective variety of dimension $m=n+2\geq 4$ and
degree $t$. Denote by $V^{(4)}$ the  general linear section of $V$
of dimension $4$. For any subvariety $X\subset V$ of dimension
$n$, denote by $X^{(h)}$ the general linear section of $X$ of
dimension $h$ ($1\leq h \leq n$), by $d$ the degree of $X$, and by
$g$ the geometric genus of $X^{(1)}$. Fix a real number $\lambda$
such that
$$
0<\lambda < \frac{1}{2\sqrt{2t}}. \tag 2.1.1
$$
Assume that $V^{(4)}$ is smooth, and that its N\acuteaccent
eron-Severi group has rank $1$. Then the set of irreducible,
projective subvarieties $X$ of $V$ of dimension $n$ such that
$X^{(2)}$ is smooth and
$$
g\leq \lambda d^\frac{3}{2},
$$
is bounded.
\endproclaim \bigskip

\demo {Proof of Theorem 2.1} We may assume $n=2$. Let $X\subset V$
be any smooth projective surface such that $g\leq \lambda
d^\frac{3}{2}$. By [CD2] we may  assume that $X$ is of general
type. Therefore  $\chi (\Cal O_X))\geq 0$, and from $(1.5)$ we get
$$
d^2/t-2K^2_X+O(d^\frac{3}{2})\leq 0. \tag 2.1.2
$$
Now consider the orthogonal decomposition $K_X=D+eH_X$ in
$NS(X)\otimes {\bold Q}$, with $D\cdot H_X=0$ and $D^2\leq 0$.
Since $e=(2g-2-d)/d$, $g\leq \lambda d^\frac{3}{2}$ and we may
assume $d > O(1)$, we have
$$
|e|\leq 2\lambda d^\frac{1}{2}.
$$
It follows that
$$
K^2_X=(D+eH_X)^2\leq e^2d\leq 4\lambda^2d^2.
$$
Using $(2.1.2)$ we have
$$
\left(\frac{1}{t}-8\lambda^2\right)d^2+O(d^\frac{3}{2})\leq 0,
$$
which, taking into account $(2.1.1)$,   proves the boundedness of
$d$. $\square$
\bigskip
\bigskip

We are in position to prove  Theorem 0.1 under the hypotheses (B),
with $a=1$ . To this aim, fix any real number $\lambda$ such that
$$
0<\lambda < \frac {1}{n!\sqrt{(2t)^n}}. \tag 2.2
$$
Notice that, taking into account Castelnuovo-Harris  bound for the
geometric genus of a projective variety [H], it suffices to prove
that the set of smooth, projective subvarieties $X$ of $V$ of
dimension $n$ with
$$
h^0(X,\Cal O_X(K_X-bH_X)) \leq \lambda d^{\frac{n}{2}}, \tag 2.3
$$
is bounded. Using our hypotheses on the homology of the ambient
variety $V$, we know that $X=\alpha H^2_V$ in $H_{2n}(V,{\bold
C})$, for some $\alpha \in {\bold Q}$. Hence we can consider the
following natural commutative diagram
$$
\matrix
 &  &  A^1(X)  \\
& & \\
 & j\swarrow &         \downarrow{\alpha H^2_X} \\
& & \\
A^{3}(V) &          @>\rho >>    &     A^{3}(X)
\endmatrix \tag 2.4
$$
\noindent where $A^{3}(V) \subset H_{2n-2}(V,{\bold C})$ denotes
the space of algebraic  classes of codimension $3$ of $V$, and
similarly for $X$, $\rho$ is the natural restriction map, and $j$
is the map which takes a cycle on $X$ and considers it as a cycle
on $V$. Notice that, by our hypotheses, $A^{3}(V)\simeq \bold C$.
Moreover the vertical map is injective by Hard Lefschetz: this is
where we use $n\geq 4$. It follows that $j$ is an isomorphism,
i.e. $A^1(X)\simeq NS(X)\otimes \bold C$ has dimension $1$. When
$a=1$, we assume $n\geq (r+1)/2$ and so we may apply  Larsen
Theorem [L]. With the same argument developed in [A], Proposition
8, pg. 69, we deduce that the Picard group of $X$ is generated by
the hyperplane section. It follows that
$$
K_X=eH_X
$$
in $Pic(X)$, with $e=(2g-2-(n-1)d)/d$. Therefore, by $(2.3)$ we
have
$$
h^0(X,\Cal O_X(K_X-bH_X))=h^0(X,\Cal O_X(e-b)) \leq \lambda
d^{\frac{n}{2}}. \tag 2.5
$$
On the other hand, by Castelnuovo Theory (see Proposition (3.23),
pg. 117 in [EH]) we know that
$$
h^0(X,\Cal O_X(e-b))\geq (e-b)^n/n! \tag 2.6
$$
(by Theorem 2.1 we may assume $e-b>0$). Combining $(2.5)$ with
$(2.6)$ we obtain
$$
g\leq \frac{(n!\lambda)^{\frac{1}{n}}}{2} d^\frac{3}{2}+O(d), \tag
2.7
$$
and the boundedness of $d$ follows by $(2.2)$ and Theorem 2.1.
This concludes the proof of  Theorem 0.1 under the hypotheses (B),
for $a=1$. $\square$
\bigskip

Now we are going to prove Theorem 0.1 under the hypotheses (B),
with $a\geq 2$. We need the following preliminary result, which
relies on Kawamata-Viehweg Vanishing Theorem, and   allows us  to
 bound $h^0(X,\Cal O_X(aK_X-bH_X))$ from below  only assuming that the
N\'eron-Severi group of $X$ has rank $1$ (compare with $(2.5)$ and
$(2.6)$ above).

\proclaim {Proposition 2.8} Let $X\subset {\bold P^r}$  be a
projective, irreducible and smooth  variety of general type, of
dimension $n\geq 2$ and degree $d$. Let $D$ be a big and nef
divisor on $V$. Assume that for suitable integers $\sigma$ and
$\tau$, the general hyperplane section $H_X$ is numerically
equivalent to $\tau D$, and the canonical divisor  $K_X$ is
numerically equivalent to $\sigma D$. Fix an integer $l>\sigma$
and define $\alpha$ and $\beta$ by dividing $l-\sigma=\alpha
\tau+\beta$, $0\leq\beta <\tau$. Put
$$
\gamma(n,d,\sigma,\tau,l) = (n-2)(\alpha
-n+1)\left[\alpha\tau^2(n-1)+n\tau(\sigma -
2l)\right]+n(n-1)l(l-\sigma-\tau(n-2)).
$$
Then one has
$$
h^0(X,\Cal O_X(lD))\geq d{{\alpha
-1}\choose{n-2}}\gamma(n,d,\sigma,\tau,l) /2\tau^2n(n-1). \tag
2.8.1
$$
\endproclaim \medskip

\noindent {\bf Remark 2.9.} (i) In $(2.8.1)$ we assume ${{\alpha
-1}\choose{n-2}}=1$ for $n=2$, and ${{\alpha -1}\choose{n-2}}=0$
for $n>2$ and $\alpha <n-1$.

(ii) Proposition 2.8 should be compared with
Koll\'ar-Luo-Matsusaka estimate [K], pg. 302, Theorem 2.15.9,
which, in our context, is not as good as inequality $(2.8.1)$.
$\square$

\demo {Proof of Proposition 2.8}  First we examine the case $n=2$.
Since $K_X\equiv \sigma D$, and $l>\sigma$, by Kawamata-Viehweg
Vanishing Theorem we have $h^i(X,\Cal O_X(lD))=0$ for any $i>0$.
Hence by Riemann-Roch Theorem  we get $h^0(X,\Cal O_X(lD))=
\chi(\Cal O_X)+l(l-\sigma)D^2/2$. Since $X$ is of general type,
then $\chi(\Cal O_X)\geq 0$ and so we deduce
$$
h^0(X,\Cal O_X(lD))\geq l(l-\sigma)D^2/2=
d\gamma(2,d,\sigma,\tau,l) /4\tau^2. \tag 2.8.2
$$
This proves Proposition 2.8 for $n=2$.

Now we are going to examine the case  $n>2$. We may  assume
$\alpha \geq n-1$. By Kawamata-Viehweg Vanishing Theorem  we have,
for any $0\leq j_1 \leq \alpha -2$,
$$
h^1(X,\Cal O_X(lD-j_1H_X))=0
$$
and
$$
h^0(X^{(n-1)},\Cal O_X(lD-j_1H_X)\otimes \Cal
O_{X^{(n-1)}})=h^0(X^{(n-1)},\Cal
O_{X^{(n-1)}}((l-j_1\tau)D{^{(n-1)}})),
$$
where $X^{(n-1)}=H_X$ denotes the general hyperplane section of
$X$ and $D^{(n-1)}$ denotes the restriction of $D$ to $X^{(n-1)}$.
Hence from the hyperplane section exact sequence
$$
0\to\Cal O_X(lD-(j_1+1)H_X)\to \Cal O_X(lD-j_1H_X)\to \Cal
O_X(lD-j_1H_X)\otimes \Cal O_{X^{(n-1)}} \to 0
$$
we get
$$
h^0(X,\Cal O_X(lD))\geq  \sum_{j_1=0}^{\alpha-2}h^0(X^{(n-1)},\Cal
O_{X^{(n-1)}}((l-j_1\tau)D{^{(n-1)}})).
$$
Iterating the previous argument for the successive linear sections
of $X$, we get
$$\multline
h^0(X,\Cal O_X(lD))\geq\\
\sum_{j_1=0}^{\alpha -2}\sum_{j_2=0}^{\alpha -3
-j_1}\dots\sum_{j_{n-2}=0}^{\alpha -(n-1)-j_1-j_2-\dots-j_{n-3}}
h^0(X^{(2)},\Cal
O_{X^{(2)}}((l-(j_1+j_2+\dots+j_{n-2})\tau)D{^{(2)}})),
\endmultline \tag 2.8.3
$$
where $X^{(2)}$ denotes the general linear section of $X$ of
dimension $2$, and $D{^{(2)}}$ the restriction of $D$ on
$X^{(2)}$. Notice that $X^{(2)}$ is of general type because $X$
is. Hence from $(2.8.2)$ we have
$$\multline
h^0(X^{(2)},\Cal
O_{X^{(2)}}((l-(j_1+j_2+\dots+j_{n-2})\tau)D{^{(2)}}))\\
\geq d[l-(j_1+j_2+\dots+j_{n-2})\tau]
[l-(j_1+j_2+\dots+j_{n-2})\tau-\sigma-(n-2)\tau)]/2\tau^2.
\endmultline \tag 2.8.4
$$
Combining $(2.8.3)$ with $(2.8.4)$ we obtain $(2.8.1)$. $\quad$
$\square$
\bigskip

We are in position to prove Theorem 0.1 under the hypotheses (B),
with $a\geq 2$. Fix any real number $\lambda$ such that
$$
0<\lambda < \frac {(a-1)^{n-1}(2a+n-2)}{2n!\sqrt{(2t)^n}}. \tag
2.10
$$
As in the case $a=1$, it suffices to prove that the set of smooth,
projective subvarieties $X$ of $V$ of dimension $n$ with
$$
h^0(X,\Cal O_X(aK_X-bH_X)) \leq \lambda d^{\frac{n+2}{2}}, \tag
2.11
$$
is bounded. As in the case $a=1$, one sees that $NS(X)$ has rank
1. Therefore, for some ample divisor $D$ on $X$, and suitable
integers $\sigma$ and $\tau$, we have $H_X\equiv \tau D$, and
$K_X\equiv \sigma D$. Put $e={\sigma}/{\tau}$ and notice that
$e=(2g-2-(n-1)d)/d$. By Theorem 2.1 we may assume $e>O(1)$, in
particular $X$ is of general type. Moreover, if we put
$l=a\sigma-b\tau$ then  we have $l>\sigma$ because $a\geq 2$.
Also, by Kawamata-Viehweg Vanishing Theorem, we have $h^0(X,\Cal
O_X(aK_X-bH_X))=h^0(X,\Cal O_X(lD))$. Hence we can apply
Proposition 2.8 and from $(2.8.1)$ we get
$$
h^0(X,\Cal O_X(aK_X-bH_X))\geq
\frac{d}{2n!}\left[(a-1)^{n-1}(2a+n-2)e^n+O(e^{n-1})\right].\tag
2.12
$$
Now fix any real number $\mu$ such that $\lambda<\mu<\frac
{(a-1)^{n-1}(2a+n-2)}{2n!\sqrt{(2t)^n}}$ (compare with $(2.10)$),
and put
$$
\lambda_1={\root{n}\of{\frac{n!\mu}{2^{n-1}(a-1)^{n-1}(2a+n-2)}}}.
\tag 2.13
$$
We have $0<\lambda_1<\frac {1}{2\sqrt{2t}}$ (see $(2.1.1)$).
Hence, from Theorem 2.1, we may assume
$g>\lambda_1d^{\frac{3}{2}}$. From $(2.12)$ we obtain
$$
h^0(X,\Cal O_X(aK_X-bH_X))\geq
\frac{(a-1)^{n-1}}{n!}2^{n-1}(2a+n-2)\lambda_1^nd^{\frac{n+2}{2}}+O(d^{\frac{n+1}{2}}).\tag
2.14
$$
In view of 2.13, comparing 2.14 with 2.11 we deduce that $d$ is
bounded.

In a similar way as under the hypotheses (A), with $n\neq 3$,
considering codimension two balanced complete intersections on $V$
(when $V$ is smooth), one proves the sharpness of Theorem 0.1
under the hypotheses (B), when $a\geq 2$, in the sense of
$(0.1.3)$.

This concludes the proof of Theorem 0.1 under the hypotheses (B),
when $a\geq 2$.
\bigskip
\bigskip

\noindent {\bf Remark 2.15.} (i) As already remarked in $(2.3)$,
we notice that, taking into account Castelnuovo-Harris  bound for
the geometric genus of a projective variety [H], in order to prove
Theorem 0.1 under the hypotheses (B), it suffices to prove the
existence of a positive real number $\lambda$ such that the set of
smooth, projective subvarieties $X$ of $V$ of dimension $n$ with $
h^0(X,\Cal O_X(aK_X-bH_X)) \leq \lambda d^{\frac{n}{2}}$, is
bounded. The corresponding claim  holds true also under the
hypotheses (D) (compare with  $(4.3.2)$ in  next Section 4). But
this is not true under the hypotheses (A) (and (C) and (E), see
below).

(ii) With obvious modification, the proof of Theorem 0.1 under the
hypotheses (B) with $a=1$ works  also for $a\geq 2$. $\square$
\bigskip
\bigskip
\bigskip

\noindent {\bf 3. The proof of Theorem 0.1 under the hypotheses
(C).}
\bigskip

As in the proof of Theorem 0.1 under the hypotheses (A), in order
to prove  Theorem 0.1 under the hypotheses (C), we may assume
$a=1$ and $b=0$. This said, the  proof consists in showing the
existence of a suitable constant $\mu>0$ such that
$$
p_g(X)-c\left(\sum_{h=1}^{n-1}p_g(X^{(h)})\right)\geq \mu
d^{(n+2)/2},
$$
for any smooth and $n-$dimensional subvariety $X\subset V$ of
degree $d>O(1)$. In order to prove this, first we compare the
geometric genus of $X$ with the arithmetic genus. This is done in
Corollary $(3.5)$ below.  To prove this result, we need some
preliminaries. First we prove Lemma $(3.1)$ and Proposition
$(3.2)$ which allow us to control the difference between the
geometric genus and the arithmetic genus in terms of the
arithmetic genera of the linear sections. Next we need Lemma
$(3.3)$ and Lemma $(3.4)$ to obtain a numerical control of these
arithmetic genera.
\bigskip
\bigskip

\proclaim {Lemma 3.1} Let $Y\subset {\bold P^r}$  be an
irreducible, smooth, projective variety of dimension $m \geq 2$.
For any $1\leq j \leq m$,  denote by $Y^{(j)}$ the general linear
section of $Y$ of dimension $j$. Then we have
$$
h^{j}(Y,\Cal O_Y)\leq p_g(Y^{(j)}).
$$
\endproclaim \bigskip

\demo {Proof} Use Hyperplane Lefschetz Theorem and induction on
$m$. $\square$
\bigskip
\bigskip

\proclaim {Proposition 3.2} Let $Y\subset {\bold P^r}$  be an
irreducible, smooth, projective variety of dimension $m \geq 1$.
Then there exist integers $c_0,\dots ,c_{m-1}$ and $d_0,\dots
,d_{m-1}$, depending only on $m$, and such that
$$
c_0+\sum_{j=1}^{m-1}c_j\chi({\Cal O_{Y^{(j)}}}) \leq
p_g(Y)+(-1)^{m+1}\chi({\Cal O_{Y}})\leq
d_0+\sum_{j=1}^{m-1}d_j\chi({\Cal O_{Y^{(j)}}}).
$$
\endproclaim \bigskip

\demo {Proof} The case $m=1$ being trivial, we may assume $m\geq
2$ and argue by induction on $m$. We have
$$\multline
p_g(Y)+(-1)^{m+1}\chi({\Cal O_{Y}})= h^m(Y, \Cal
O_Y)+(-1)^{m+1}\sum_{j=0}^{m}(-1)^jh^j(Y, \Cal O_Y)\\
=(-1)^{m+1}\sum_{j=0}^{m-1}(-1)^jh^j(Y, \Cal O_Y).
\endmultline
$$
By Lemma $(3.1)$ we deduce
$$
-(1+\sum_{j=1}^{m-1}p_g(Y^{(j)}))\leq
 p_g(Y)+(-1)^{m+1}\chi({\Cal
O_{Y}}) \leq 1+\sum_{j=1}^{m-1}p_g(Y^{(j)}),
$$
and  our claim follows by using the induction hypothesis on the
general linear sections $Y^{(1)},\dots,Y^{(m-1)}$. $\square$
\bigskip
\bigskip

\proclaim {Lemma 3.3} Let $V\subset {\bold P^r}$  be an
irreducible, smooth, projective variety of dimension $n+2 \geq 7$.
Assume that for $1\leq i\leq 3$, any algebraic class in
$H^{2i}(V,{\bold C})$ is a multiple of $H^{i}_V$, where $H_V$ is
the hyperplane section of $V$. Let $X\subset V$  be an
irreducible, smooth, projective subvariety of dimension $n$, with
degree $d$, and assume that $d>O(1)$. Then
$$
c_1^2(N_{X,V})\geq 2c_2(N_{X,V})>0.
$$
\endproclaim \bigskip

\demo {Proof} We already know that, by our hypotheses on the
cohomology of $V$, $X$ is numerically subcanonical (see $(2.4)$
before). Hence we have $K_X=eH_X$ in $H^{2}(X,{\bold C})$, where
$e=(2g-2-(n-1)d)/d$ (as usual, $g$ denotes the linear genus of
$X$). Therefore we have $c_1(N_{X,V})=n_1H_X$ in $H^{2}(X,{\bold
C})$, where
$$
n_1=e-q,
$$
with $K_V=qH_V$ in $H^{2}(V,{\bold C})$. On the other hand, from
the self-intersection formula, we have $c_2(N_{X,V})=n_2H_X^2$ in
$H^{4}(X,{\bold C})$, where
$$
n_2=d/t
$$
($t=$ degree of $V$). Hence we only have to prove that
$$
n_1^2\geq 2n_2.
$$
To this purpose, fix a constant $\gamma >>0$ such that the twisted
tangent bundle $T_V(\gamma)$ is globally generated. Then also
$N_{X,V}(\gamma)$ is, and therefore we have the positivity of its
Segre class (see [FL])
$$\multline
-s_5(N_{X,V}(\gamma))=c_1^5(N_{X,V}(\gamma))-4c_1^3(N_{X,V}(\gamma))c_2(N_{X,V}(\gamma))\\
+3c_1(N_{X,V}(\gamma))c_2^2(N_{X,V}(\gamma))\geq 0.
\endmultline \tag 3.3.1
$$
Taking into account that $c_1(N_{X,V}(\gamma))=(n_1+2\gamma)H_X$,
$c_2(N_{X,V}(\gamma))=(n_2+\gamma n_1+\gamma^2)H_X^2$, and that by
Theorem 2.1 we may assume
$$
n_1+2\gamma>0  \quad\quad{\text {and}}\quad\quad n_2+\gamma
n_1+\gamma^2>0,\tag 3.3.2
$$
the previous inequality $(3.3.1)$ yields
$$
y^2-4y+3\geq 0,
$$
where
$$
y=(n_1+2\gamma)^2/(n_2+\gamma n_1+\gamma^2).
$$
It follows that either $y\leq 1$ or $y\geq 3$. From the positivity
$$
-s_3(N_{X,V}(\gamma))=\left[(n_1+2\gamma)^3-2(n_1+2\gamma)(n_2+\gamma
n_1+\gamma^2)\right]\cdot H_X^3\geq 0,
$$
we have  $2(n_2+\gamma n_1+\gamma^2)\leq (n_1+2\gamma)^2$. Hence,
if $y\leq 1$ then $n_2+\gamma n_1+\gamma^2\leq 0$, which is in
contrast with $(3.3.2)$. We deduce that $y\geq 3$, i.e.
$$
(3n_2-n_1^2)-\gamma n_1-\gamma^2\leq 0.
$$
Now suppose $n_1^2\leq 2n_2$. From the previous inequality we have
$$
n_1^2/2-\gamma n_1-\gamma^2\leq 0.
$$
We obtain $n_1\leq O(1)$ which, by Theorem 2.1, implies that $d$
is bounded. This is in contrast with our hypotheses $d>O(1)$.
Hence we must have $n_1^2\geq 2n_2$,  and this concludes the proof
of  Lemma $(3.3)$. $\square$
\bigskip
\bigskip

\proclaim {Lemma 3.4} Let $V\subset {\bold P^r}$  be an
irreducible, smooth, projective variety of dimension $n+2 \geq 7$.
Assume that, for $1\leq i\leq 3$, any algebraic class in
$H^{2i}(V,{\bold C})$ is a multiple of $H^{i}_V$, where $H_V$ is
the hyperplane section of $V$. Moreover assume that, for $1\leq
i\leq n$, there exist rational numbers $v_1,\dots,v_n$ such that
$c_i(T_V)=v_iH_V^i$ in $H^{2i}(V,{\bold C})$. For any irreducible,
smooth, projective subvariety $X\subset V$ of dimension $n$, with
degree $d$ and linear genus $g$, put $n_1=(2g-2-(n-1)d)/d +v_1$.
Assume $d>O(1)$. Then $n_1>0$ and, for any $1\leq h\leq n$, one
has
$$
|\chi({\Cal O_{X^{(h)}}})|\leq dO(n_1^h).
$$
\endproclaim \bigskip

\demo {Proof} Fix an integer $h\in \{1,\dots,n\}$. By
Hirzebruch-Riemann-Roch Theorem we know that $\chi({\Cal
O_{X^{(h)}}})$ is equal to the degree of the top Todd class of
$X^{(h)}$. This is the degree of a certain linear combination,
with coefficients depending only on $h$, of monomials like
$$
c_1^{\alpha_1}(T_{X^{(h)}})\cdot c_2^{\alpha_2}(T_{X^{(h)}})\dots
c_h^{\alpha_h}(T_{X^{(h)}}),
$$
with $\alpha_1,\dots,\alpha_h$ non-negative integers, and
$h=\sum_{j=1}^h j\alpha_j$. From  the natural exact sequence
$$
0\to T_{X^{(h)}}\to T_{V^{(h+2)}}\otimes \Cal O_{X^{(h)}}\to
N_{{X^{(h)}},{V^{(h+2)}}}\to 0,
$$
we may compute the Chern classes of ${X^{(h)}}$ in terms of the
Chern classes of $T_{V^{(h+2)}}\otimes \Cal O_{X^{(h)}}$, and the
Chern classes of $N_{{X^{(h)}},{V^{(h+2)}}}$. With the same
notation as in the proof of Lemma $(3.3)$, we have
$$
c_1(N_{{X^{(h)}},{V^{(h+2)}}})=n_1H_{X^{(h)}}\quad {\text{and}}
\quad c_2(N_{{X^{(h)}},{V^{(h+2)}}})=n_2H_{X^{(h)}}^2.
$$
Taking into account our hypotheses on the Chern classes of $V$, it
follows that  $\chi({\Cal O_{X^{(h)}}})$ is equal to a certain
linear combination, with coefficients depending only on $h$, of
monomials like
$$
dn_1^{\alpha_1}n_2^{\alpha_2}v_1^{\beta_1}v_2^{\beta_2}\dots
v_h^{\beta_h},
$$
with $\alpha_1,\alpha_2,\beta_1,\beta_2,\dots,\beta_h$
non-negative integers such that $\alpha_1+2\alpha_2+\sum_{j=1}^h
j\beta_j\leq h$. From Lemma $(3.3)$ we know that if $d>O(1)$ then
$$
n_1^2\geq 2n_2 > 1.
$$
Hence, for a suitable constant $K$, we have
$$
|dn_1^{\alpha_1}n_2^{\alpha_2}v_1^{\beta_1}v_2^{\beta_2}\dots
v_h^{\beta_h}|\leq Kdn_1^h.
$$
This proves that $|\chi({\Cal O_{X^{(h)}}})|\leq dO(n_1^h)$.
$\square$
\bigskip
\bigskip

As a consequence of Proposition $(3.2)$, Lemma $(3.3)$ and Lemma
$(3.4)$, we obtain the following
\bigskip
\bigskip

\proclaim {Corollary 3.5} Let $V\subset {\bold P^r}$  be an
irreducible, smooth, projective variety of dimension $n+2 \geq 7$.
Assume that, for $1\leq i\leq 3$, any algebraic class in
$H^{2i}(V,{\bold C})$ is a multiple of $H^{i}_V$, where $H_V$ is
the hyperplane section of $V$. Moreover assume that, for $1\leq
i\leq n$, there exist rational numbers $v_1,\dots,v_n$ such that
$c_i(T_V)=v_iH_V^i$ in $H^{2i}(V,{\bold C})$. Fix an integer
$c\geq 0$. For any irreducible, smooth, projective subvariety
$X\subset V$ of dimension $n$, with degree $d$ and linear genus
$g$, put $n_1=(2g-2-(n-1)d)/d +v_1$. Assume $d>O(1)$. Then $n_1>0$
and
$$
p_g(X)-c\left(\sum_{h=1}^{n-1}p_g(X^{(h)})\right)\geq
(-1)^n\chi({\Cal O_{X}})+dO(n_1^{n-1}). \quad \square
$$
\endproclaim
\bigskip
\bigskip

We are in position to prove Theorem 0.1 under the hypotheses (C).
To this purpose, let $X\subset V$ be a smooth subvariety of
codimension $2$. By Corollary $(3.5)$, in order to bound
$p_g(X)-c\left(\sum_{h=1}^{n-1}p_g(X^{(h)})\right)$ from below, it
is enough to  bound $(-1)^n\chi({\Cal O_{X}})$ from below. We keep
the notation we introduced in the proof of Lemma $(3.3)$. As in
the proof of this lemma, using Hirzebruch-Riemann-Roch Theorem,
one may write
$$
(-1)^n\chi({\Cal O_{X}})=(-1)^nd\left(U+R\right),
$$
where $U$ is a linear combination, with coefficients depending
only on $n$, of monomials like
$$
n_1^{\alpha_1}n_2^{\alpha_2},
$$
with $\alpha_1,\alpha_2$ non-negative integers such that
$\alpha_1+2\alpha_2=n$, and $R$ is a linear combination, with
coefficients depending only on $n$, of monomials like
$$
n_1^{\alpha_1}n_2^{\alpha_2}v_1^{\beta_1}v_2^{\beta_2}\dots
v_h^{\beta_h},
$$
with $\alpha_1,\alpha_2,\beta_1,\beta_2,\dots,\beta_h$
non-negative integers such that $\alpha_1+2\alpha_2+\sum_{j=1}^h
j\beta_j=n$ and $\alpha_1+2\alpha_2<n$. As in the proof of Lemma
$(3.3)$ one sees that
$$
(-1)^n dR=dO(n_1^{n-1}).
$$
Summing up we obtain
$$
p_g(X)-c\left(\sum_{h=1}^{n-1}p_g(X^{(h)})\right)\geq (-1)^ndU+
dO(n_1^{n-1}).\tag 3.6
$$
Now define $\nu$ and $\rho$ by dividing
$$
n=2\nu+\rho, \quad 0\leq \rho \leq 1,
$$
and denote by
$$
x_0,x_1,\dots,x_{\nu}
$$
the integers, depending only on $n$, such that
$$
(-1)^nU={\frac{1}{(n+2)!}}\left(\sum_{j=0}^\nu
x_jn_1^{n-2j}n_2^j\right). \tag 3.7
$$
Now we need the following numerical lemma. We will prove it later.
\bigskip \bigskip

\proclaim {Lemma  3.8} With the same notation as before, denote by
$$
q(y)=\sum_{j=0}^\nu x_jy^{\nu-j},
$$
and by $q^{(j)}(y)$ its $j-$th derivative. Then one has
$q^{(j)}(2)>0$, for any $0\leq j\leq \nu$. $\square$
\endproclaim
\bigskip
\bigskip
As a consequence we have the following
\bigskip \bigskip

\proclaim {Corollary 3.9} There exists a rational number $\epsilon
>0$, depending only on $n$, such that, if we put
$$
q_{\epsilon}(y)=q(y)-\epsilon y^{\nu},
$$
then one has $q_{\epsilon}^{(j)}(2)>0$, for any $0\leq j\leq \nu$.
In particular, if $y\geq 2$ then $q_{\epsilon}(y)\geq
q_{\epsilon}(2)$.$\square$
\endproclaim
\bigskip
\bigskip

Now, continuing our computation,  from $(3.6)$ and $(3.7)$, and
taking into account Lemma $(3.3)$, Lemma $(3.8)$ and Corollary
$(3.9)$, we may write
$$\multline
p_g(X)-c\left(\sum_{h=1}^{n-1}p_g(X^{(h)})\right)\geq
{\frac{d}{(n+2)!}}n_1^{\rho}n_2^{\nu}\left[q_{\epsilon}\left({\frac{n_1^2}{n_2}}\right)
+\epsilon\left({\frac{n_1^2}{n_2}}\right)^{\nu}\right]+dO(n_1^{n-1})\\
={\frac{d}{(n+2)!}}n_1^{\rho}n_2^{\nu}q_{\epsilon}\left({\frac{n_1^2}{n_2}}\right)+
{\frac{d}{(n+2)!}}\epsilon n_1^n+dO(n_1^{n-1})\\
={\frac{d}{(n+2)!}}n_1^{\rho}n_2^{\nu}q_{\epsilon}\left({\frac{n_1^2}{n_2}}\right)+
{\frac{d}{(n+2)!}}n_1^{n-1}\left[\epsilon n_1+O(1)\right] \geq
{\frac{d}{(n+2)!}}n_1^{\rho}n_2^{\nu}q_{\epsilon}\left({\frac{n_1^2}{n_2}}\right)\\
\geq {\frac{d}{(n+2)!}}n_1^{\rho}n_2^{\nu}q_{\epsilon}(2) \geq
{\frac{d}{(n+2)!}}n_2^{\nu+\rho/2}q_{\epsilon}(2)
={\frac{q_{\epsilon}(2)}{(n+2)!\sqrt{t^n}}}d^{{\frac{n}{2}}+1}.
\endmultline
$$
Since ${\frac{q_{\epsilon}(2)}{(n+2)!\sqrt{t^n}}}$ is a constant
$>0$, the previous inequality proves Theorem 0.1 under the
hypotheses (C) (as before, considering codimension two balanced
complete intersections on $V$, one proves the sharpness in the
sense of $(0.1.3)$).

\smallskip
It remains  to prove the numerical Lemma $(3.8)$. To this aim,
keep all the notation we introduced before. Since the polynomial
$q(y)$ depends only on $n$, in order to compute it, we may use a
complete intersection $X\subset V$ of type $(u,v)$. In this case
we have
$$
n_1=u+v \quad {\text{and}} \quad n_2=uv.
$$
Moreover
$$
\chi({\Cal O_{X}})=p_V(0)-p_V(-u)-p_V(-v)+p_V(-u-v)
$$
($p_V(l)$ = Hilbert polynomial of $V$). We may compute the term
$(-1)^ndU$ assuming that all Chern classes of $V$ are $0$. By
Hirzebruch-Riemann-Roch Theorem again, this implies that we may
assume
$$
p_V(l)=t{\frac{l^{n+2}}{(n+2)!}}
$$
($t$ = degree of $V$). Since $d=uvt$ we obtain
$$
(-1)^ndU={\frac{d}{(n+2)!}}\left[{\frac{(u+v)^{n+2}-u^{n+2}-v^{n+2}}{uv}}\right].
$$
In other words, the coefficients $x_0,x_1,\dots,x_{\nu}$ of the
polynomial $q(y)$ are defined by the identity
$$
{\frac{(u+v)^{n+2}-u^{n+2}-v^{n+2}}{uv}}=\sum_{j=0}^{\nu}x_j(u+v)^{n-2j}(uv)^j.
\tag 3.8.1
$$
Now, if we put
$$
y={\frac{(u+v)^2}{uv}},
$$
we have
$$
q(y)=y^{\nu +1}-{\frac{u^{n+2}+v^{n+2}}{(u+v)^{\rho}(uv)^{\nu
+1}}}. \tag 3.8.2
$$
Fix a real number $3/2\leq y\leq 5/2$. Then we have
$y={\frac{(u+v)^2}{uv}}$ with
$u={\frac{\left[y-2+\sqrt{-1}\sqrt{y(4-y)}\right]}{2}}$ and $v=1$.
Notice that since $|u|=1$, then we have
$$
u=exp({\sqrt{-1}\theta})
$$
for some real number $\theta$. In particular we have
$$
y-2=2\cos\theta  \quad{\text{and}}\quad \sqrt{y(4-y)}=2\sin\theta.
\tag 3.8.3
$$
From $(3.8.2)$ we get
$$
q(y)=y^{\nu +1}-{\frac{u^{n+2}+1}{(u+1)^{\rho}u^{\nu +1}}}, \tag
3.8.4
$$
and so, taking into account that if $y=2$ then $u=\sqrt{-1}$, we
deduce
$$
q(2) = \cases
2^{\nu +1} & {\text{if $n\equiv 0\quad mod(4)$}},\\
2^{\nu +1}+(-1)^{{\frac{n-1}{4}}} & {\text{if $n\equiv 1\quad mod(4)$}},\\
2^{\nu +1}+2(-1)^{{\frac{n+6}{4}}} & {\text{if $n\equiv 2\quad mod(4)$}},\\
2^{\nu +1}+(-1)^{{\frac{n+5}{4}}} & {\text{if $n\equiv 3\quad
mod(4)$}}.
\endcases
$$
This proves that $q(2)=q^{(0)}(2)>0$.

Moreover, from $(3.8.1)$, one sees that $x_0=n+2$, and therefore
$q^{(\nu)}(2)>0$.

It remains to evaluate the derivatives $q^{(j)}(2)$, for $1\leq
j\leq \nu -1$.

First we examine the case $n$ is even. Hence $n=2\nu$ and
$\rho=0$. In this case, from $(3.8.4)$, we may write
$$
q(y)=y^{\nu +1}-\left(u^{\nu+1}+{{\frac{1}{u^{\nu+1}}}}\right)=
y^{\nu
+1}-\left[exp(\sqrt{-1}(\nu+1)\theta)+exp(-\sqrt{-1}(\nu+1)\theta)\right]
$$
and so we get
$$
q(y)=y^{\nu +1}-2\cos(\nu+1)\theta. \tag 3.8.5
$$
One may write  $\cos(\nu+1)\theta$ as a polynomial in $\cos\theta$
of degree $\nu+1$, i.e.
$$
\cos(\nu+1)\theta=\sum_{l=0}^{\nu+1}\beta_l(\cos\theta)^l,
$$
with suitable integer numbers $\beta_0,\dots,\beta_{\nu+1}$. By
$(3.8.3)$ we obtain the following polynomial identity
$$
q(y)=y^{\nu
+1}-2\left[\sum_{l=0}^{\nu+1}\beta_l{{\frac{(y-2)^l}{2^l}}}\right].
$$
From this formula we deduce, for $1\leq j \leq \nu -1$,
$$
q^{(j)}(2)={{\frac{j!}{2^{j-1}}}}\left[2^{\nu}{{\nu+1}\choose{j}}-\beta_j\right].
\tag 3.8.6
$$
Now we need the following
\bigskip \bigskip

\proclaim {Sublemma} With the same notation as before, assume
$\nu\geq 0$. Then for any $0\leq l \leq \nu+1$ one has
$$
|\beta_l|\leq 2^{l}{{\nu+1}\choose{l}}.
$$
\endproclaim
\bigskip

\demo {Proof of the Sublemma} The case $0\leq \nu \leq 1$ being
trivial, we may assume $\nu \geq 2$ and argue by induction on
$\nu$. From the identity
$$
\cos(\nu+1)\theta=2\cos\theta\cos \nu\theta-\cos(\nu-1)\theta \tag
3.8.7
$$
and the induction hypotheses, we deduce that
$$
\beta_{\nu+1}=2^{\nu}, \quad \beta_{\nu}=0,\quad {\text{and}}
\quad |\beta_0|\leq 1.
$$
Hence we only have to estimate $|\beta_l|$ for  $1\leq l\leq
\nu-1$. To this purpose, put
$$
\cos\nu\theta=\sum_{l=0}^{\nu}\gamma_l(\cos\theta)^l \quad
{\text{and}} \quad
\cos(\nu-1)\theta=\sum_{l=0}^{\nu-1}\delta_l(\cos\theta)^l.
$$
From $(3.8.7)$ we have $\beta_l=2\gamma_{l-1}-\delta_l$.
Therefore, using the induction hypotheses, we get
$$
|\beta_l|\leq 2|\gamma_{l-1}|+|\delta_l|\leq
2^{l}{{\nu}\choose{l-1}}+2^{l}{{\nu-1}\choose{l}}\leq
2^{l}{{\nu+1}\choose{l}}.\quad \square
$$
\bigskip
\bigskip
\noindent Continuing the computation of $q^{(j)}(2)$ from
$(3.8.6)$ and using the sublemma, we get
$$
q^{(j)}(2)\geq 2j!{{\nu+1}\choose{j}}(2^{\nu -j}-1)>0.
$$
This concludes the proof of Lemma $(3.8)$ in the case $n$ even.

Finally assume $n$ is odd, hence $n=2\nu+1$ and $\rho=1$. In this
case, from $(3.8.4)$, we have
$$
q(y)=y^{\nu+1}-(-1)^{\nu+1}-\left[\sum_{l=1}^{\nu+1}(-1)^{\nu+l-1}(u^l+u^{-l})\right].
$$
Therefore we get
$$
q(y)=y^{\nu+1}+(-1)^{\nu+2}-2\left[\sum_{l=1}^{\nu+1}(-1)^{\nu+l-1}\cos
l\theta\right]. \tag 3.8.8
$$
A direct computation proves Lemma $(3.8)$ for $2\leq \nu \leq 4$.
Hence we may  assume $\nu >4$ and argue by induction on $\nu$.
From $(3.8.8)$ we may write
$$
q(y)=r(y)+\left\{y^{\nu-1}+(-1)^{\nu}-2\left[\sum_{l=1}^{\nu-1}(-1)^{\nu+l-3}\cos
l\theta\right]\right\},
$$
where
$$
r(y)=\left[y^{\nu+1}-2\cos(\nu+1)\theta\right]-\left[y^{\nu-1}-2\cos(\nu\theta)\right].
$$
By induction hypotheses, all the derivatives
$q^{(j)}(2)-r^{(j)}(2)$ are $\geq 0$ for any $j\geq 0$. Therefore
we only have to prove that $r^{(j)}(2)>0$ for $1\leq j\leq \nu-1$.
This follows by rewriting $r(y)$ as
$$
r(y)=\left[y^{\nu+1}-2\cos(\nu+1)\theta\right]-\left[y^{\nu}-2\cos(\nu\theta)\right]+y^{\nu}-y^{\nu-1},
$$
and using a similar computation as in the case $n$ even (compare
with $(3.8.5)$). This concludes the proof of Lemma $(3.8)$.
\bigskip

\noindent {\bf Remark 3.9.} From $(3.8.1)$ one obtains explicit
formulae for $x_0,x_1,\dots,x_{\nu}$, i.e. one has, for $0\leq j
\leq \nu$,
$$
x_j=(-1)^j{{n-j+2}\choose{j+1}}+\sum_{l=0}^{j-1}(-1)^{j+1+l}(l+1){{n-j+2}\choose{j-1-l}}.
$$
In particular we see that
$$
x_0=n+2,
$$
and deduce $q(y)$ for low $n$. For example we have
$$
q(y) = \cases
7y^2-14y+7 & {\text{if $n=5$}},\\
8y^3-20y^2+16y-2 & {\text{if $n=6$}},\\
9y^3-27y^2+30y-9 & {\text{if $n=7$}},\\
10y^4-35y^3+50y^2-25y+2 & {\text{if $n=8$}}.
\endcases
$$
\bigskip
\bigskip
\bigskip

\noindent {\bf 4. The proof of  Theorem 0.1 under the hypotheses
(D).}
\bigskip

First we prove a boundedness result for subcanonical subvarieties
(see Theorem 4.1 and Corollary 4.3 below), which does not need the
assumption $n\geq \frac{r+1}{2}$ appearing in the hypotheses (D).

\medskip

\proclaim {Theorem 4.1} Let $V\subset {\bold P^r}$  be an
irreducible, projective variety of dimension $m$ and degree $t$.
Let $n$ be an integer with $n<m\leq 2n$, and put $k=m-n$. Assume
that  for some $1\leq h \leq n$ with $h\geq k$, the general linear
section $V^{(h+k)}$ of dimension $h+k$ of $V$ is smooth, and that
either any algebraic class in $H^{2i}(V^{(h+k)},{\bold C})$ is a
multiple of $H^i_{V^{(h+k)}}$, for $i\in \lbrace 1,\, k \rbrace$,
or $k$ is even and $H^{2i}(V^{(h+k)},{\bold C})\simeq \bold C$ for
any $i=1,\dots,k-1$. For any subvariety $X\subset V$ of dimension
$n$ and any $1\leq h \leq n$, denote by $X^{(h)}$  the general
$h$-dimensional linear section of $X$. Put $d=deg(X)$,
$g=p_g(X^{(1)})$,  and fix a real number $\lambda$ such that
$$
0<\lambda < \frac{1}{2\root{k}\of{t}}. \tag 4.1.1
$$
\medskip
Then the set of irreducible, projective subvarieties $X$ of $V$ of
dimension $n$ such that $X^{(h)}$ is smooth, numerically
subcanonical, and such that
$$
g\leq \lambda d^\frac{k+1}{k},
$$
is bounded.
\endproclaim \bigskip

\demo {Proof of Theorem 4.1}  We may assume $n=h$ by taking
$V=V^{(h+k)}$. Therefore $V$ is smooth. Let $X\subset V$ be any
smooth projective subvariety of dimension $n$, of degree $d$,
numerically subcanonical, and such that $g\leq \lambda
d^\frac{k+1}{k}$. We have
$$
K_X=eH_X
$$
in $H^{2}(X,{\bold Q})$, with $e=(2g-2-(n-1)d)/d$. Since we may
assume $d>O(1)$, we have
$$
|e|\leq 2\lambda d^\frac{1}{k}. \tag 4.1.2
$$
Now let $q$ and $\gamma$ be rational numbers such that $K_V=qH_V$
in $H^{2}(V,{\bold C})$ and $T_V(\gamma)$ is globally generated.
Then also $N_{X,V}(\gamma)$ is globally generated. Since
$$
c_1(N_{X,V}(\gamma))=(k\gamma +e-q)H_X, \tag 4.1.3
$$
by [FL]  we deduce $0 \leq k\gamma +e-q$. Notice that taking
$\gamma >O(1)$ we may assume
$$
1\leq k\gamma +e-q. \tag 4.1.4
$$
We need the following lemma. We will prove it later.

\medskip
\proclaim {Lemma 4.2} With the same assumption as before, for any
$i=1,\dots,k$ one has
$$
\multline d(k\gamma+e-q)^{n-k}O(d^{\frac{i-1}{k}}) \leq (\gamma
H_X)^{k-i}c_i(N_{X,V})c_1(N_{X,V}(\gamma))^{n-k}\\
 \leq
d(k\gamma+e-q)^{n-k}\left[\gamma^{k-i}e^i+O(d^{\frac{i-1}{k}})\right].\\
\endmultline \tag 4.2.1
$$
\endproclaim \medskip

\noindent Using the previous Lemma 4.2 for $i=k$, we have
$$
c_k(N_{X,V})c_1(N_{X,V}(\gamma))^{n-k}
 \leq
d(k\gamma+e-q)^{n-k}\left[e^k+O(d^\frac{k-1}{k})\right]. \tag
4.1.5
$$
On the other hand, if any algebraic class in $H^{2k}(V,{\bold C})$
is a multiple of $H^k_{V}$, then $X=(d/t)H^k_V$ in
$H^{2k}(V,{\bold C})$, and therefore by the self-intersection
formula we get
$$
c_k(N_{X,V})=(d/t)H^k_X.
$$
By (4.1.3) we deduce
$$
c_k(N_{X,V})c_1(N_{X,V}(\gamma))^{n-k}=d^2(k\gamma+e-q)^{n-k}/t.
\tag 4.1.6
$$
If $k$ is even and $H^{2i}(V,{\bold C})\simeq \bold C$ for any
$i=1,\dots,k-1$,  Lefschetz Hyperplane Theorem and  Hodge-Riemann
bilinear relations imply that the intersection form on
$H^{2k}(V^{(2k)},{\bold R})$ is positive definite. Hence we have
$$
(X^{(k)}-(d/t)H^k_{V^{(2k)}})^2\geq 0.
$$
Using the self-intersection formula again we get
$$
c_k(N_{X^{(k)},V^{(k)}})=(X^{(k)})^2\geq d^2/t.
$$
Using $(4.1.3)$ it follows that
$$\multline
c_k(N_{X,V})c_1(N_{X,V}(\gamma))^{n-k}=
(k\gamma+e-q)^{n-k}c_k(N_{X,V})H^{n-k}_X \\
=(k\gamma+e-q)^{n-k}c_k(N_{X^{(k)},V^{(k)}})\geq
d^2(k\gamma+e-q)^{n-k}/t,\\
\endmultline  \tag 4.1.7
$$
which, by $(4.1.6)$, holds in {\it{any}} case.

Summing up, from $(4.1.2)$, $(4.1.4)$, $(4.1.5)$ and $(4.1.7)$ we
obtain
$$
d/t\leq e^k+O({d^\frac{k-1}{k}})\leq (2\lambda)^k d
+O({d^\frac{k-1}{k}}).
$$
Therefore we have
$$
\left[\frac{1}{t}-(2\lambda)^k\right]d+O({d^\frac{k-1}{k}})\leq 0,
$$
which, taking into account $(4.1.1)$, proves the boundedness of
$d$. This concludes the proof of Theorem 4.1. $\square$
\bigskip

Now we are going  to prove Lemma 4.2. To this purpose, we argue by
induction on $i$. When $i=1$, from $(4.1.3)$ we have
$$
(\gamma H_X)^{k-1}c_1(N_{X,V})c_1(N_{X,V}(\gamma))^{n-k}=
\gamma^{k-1}d(e-q)(k\gamma+e-q)^{n-k},
$$
which, taking into account $(4.1.4)$, proves Lemma 4.2 for $i=1$.
Assume then $2\leq i\leq k$. From the formula
$$
c_i(N_{X,V}(\gamma))=\sum_{j=0}^{i}{{k-j}\choose{i-j}}(\gamma
H_X)^{i-j}c_j(N_{X,V}),
$$
intersecting with $(\gamma H_X)^{k-i}c_1(N_{X,V}(\gamma))^{n-k}$,
we get
$$
\multline
 (\gamma
H_X)^{k-i}c_i(N_{X,V})c_1(N_{X,V}(\gamma))^{n-k}= (\gamma
H_X)^{k-i}c_i(N_{X,V}(\gamma))c_1(N_{X,V}(\gamma))^{n-k}\\
-\sum_{j=1}^{i-1}{{k-j}\choose{i-j}}(\gamma
H_X)^{k-j}c_j(N_{X,V})c_1(N_{X,V}(\gamma))^{n-k}
-{{k}\choose{i}}\gamma^kd(k\gamma+e-q)^{n-k}.\\
\endmultline
$$
Using induction and $(4.1.2)$ we have
$$
{{k-j}\choose{i-j}}(\gamma
H_X)^{k-j}c_j(N_{X,V})c_1(N_{X,V}(\gamma))^{n-k}=
d(k\gamma+e-q)^{n-k}O(d^{\frac{j}{k}})
$$
for any $1\leq j \leq i-1$. It follows that
$$
\multline
 (\gamma
H_X)^{k-i}c_i(N_{X,V})c_1(N_{X,V}(\gamma))^{n-k}\\= (\gamma
H_X)^{k-i}c_i(N_{X,V}(\gamma))c_1(N_{X,V}(\gamma))^{n-k} +
d(k\gamma+e-q)^{n-k}O(d^{\frac{i-1}{k}}).
\endmultline \tag 4.2.2
$$
Now notice that since $N_{X,V}(\gamma)$ is globally generated, by
[FL] we have
$$
0\leq (\gamma
H_X)^{k-i}c_i(N_{X,V}(\gamma))c_1(N_{X,V}(\gamma))^{n-k}.\tag
4.2.3
$$
Hence the left hand side inequality in $(4.2.1)$ holds. On the
other hand one has (use [FL] and the proof of Proposition in [S])
$$
c_i(N_{X,V}(\gamma))c_1(N_{X,V}(\gamma))^{n-i}\leq
c_1(N_{X,V}(\gamma))^{n},
$$
from which, using $(4.1.2)$, $(4.1.3)$ and $(4.1.4)$, we deduce
$$
\multline (\gamma
H_X)^{k-i}c_i(N_{X,V}(\gamma))c_1(N_{X,V}(\gamma))^{n-k}
=c_i(N_{X,V}(\gamma)){\gamma}^{k-i}(k\gamma+e-q)^{n-k}H_X^{n-i}\\
={\gamma}^{k-i}(k\gamma+e-q)^{i-k}c_i(N_{X,V}(\gamma))(k\gamma+e-q)^{n-i}H_X^{n-i}\\
={\gamma}^{k-i}(k\gamma+e-q)^{i-k}c_i(N_{X,V}(\gamma))c_1(N_{X,V}(\gamma))^{n-i}\\
\leq {\gamma}^{k-i}(k\gamma+e-q)^{i-k}c_1(N_{X,V}(\gamma))^{n}
={\gamma}^{k-i}d(k\gamma+e-q)^{n+i-k}\\
=d(k\gamma+e-q)^{n-k}\left[\gamma^{k-i}e^i+O(d^{\frac{i-1}{k}})\right].\\
\endmultline
$$
Comparing the previous inequality with $(4.2.2)$ and $(4.2.3)$, we
get $(4.2.1)$. This concludes the proof of Lemma 4.2. $\square$
\bigskip

\proclaim {Corollary 4.3} Let $V\subset {\bold P^r}$  be an
irreducible, projective variety of dimension $m$ and degree $t$.
Let $n$ be an integer with $n<m\leq 2n$, and put $k=m-n$. Assume
that  for some $1\leq h \leq n$ with $h\geq k$, the general linear
section $V^{(h+k)}$ of dimension $h+k$ of $V$ is smooth, and that
either any algebraic class in $H^{2i}(V^{(h+k)},{\bold C})$ is a
multiple of $H^i_{V^{(h+k)}}$, for $i\in \lbrace 1,\, k \rbrace$,
or $k$ is even and $H^{2i}(V^{(h+k)},{\bold C})\simeq \bold C$ for
any $i=1,\dots,k-1$.

Fix integers $a,\,b,\,c \in {\bold Z}$ with $a\geq 1$, and put
$\epsilon(a)=min \lbrace \frac{n}{k}+1,\, \frac{n}{k}+a-1\rbrace$.
For any $n$-dimensional subvariety $X$ of $V$ denote by $X^{(h)}$
($K_X$ resp.) the general linear section of dimension $h$ (the
canonical divisor resp.) of $X$. Denote by $p_g(X^{(h)})$ the
geometric genus of $X^{(h)}$. Put $H_X= X^{(n-1)}$ and $d=deg(X)$.

Then there exists a strictly positive real number $\lambda
>0$, depending only on $n$, $a$ and the ambient variety $V$, such
that the set of irreducible, smooth, subcanonical, projective
subvarieties $X$ of $V$ of dimension $n$ such that
$$
h^0(X,\Cal O_X(aK_X-bH_X)) \leq \lambda
d^{\epsilon(a)}+c\left(\sum_{1\leq h <
{\epsilon(a)}-1}p_g(X^{(h)})\right),
$$
is bounded.
\endproclaim \bigskip

\demo {Proof of Corollary 4.3} First we analyze the case $a=1$.
Fix any real number $\lambda$ such that
$$
0<\lambda < \frac {1}{n!\root{k}\of{t^n}}, \tag 4.3.1
$$
where $t=deg(V)$. As in the proof of Theorem 0.1 under the
hypotheses (B) (see Section 2), one sees that in order to prove
the claim, it suffices to prove that the set of smooth,
projective, subcanonical  subvarieties $X$ of $V$ of dimension $n$
with
$$
h^0(X,\Cal O_X(K_X-bH_X)) \leq \lambda d^{\frac{n}{k}}, \tag 4.3.2
$$
is bounded. Using a similar argument as in the proof of $(2.7)$,
one sees that for such subvarieties $X$ one has
$$
g\leq \frac{(n!\lambda)^{\frac{1}{n}}}{2} d^\frac{k+1}{k}+O(d).
$$
So the boundedness of $d$ follows by $(4.3.1)$ and Theorem 4.1.

The case $a\geq 2$ follows using Proposition 2.8 in a similar
manner as in the proof of Theorem 0.1 under the hypotheses (B),
with $a\geq 2$. $\square$
\bigskip
\bigskip

We are in position to prove Theorem 0.1 in the hypotheses (D). Fix
any smooth subvariety $X$ of $V$ of dimension $n$. Using  our
hypotheses on the homology of the ambient variety $V$,  the
inequality $n\geq (m+2)/2$, and a Barth-Lefschetz type of argument
(see diagram $(2.4)$), one sees that $NS(X)\otimes \bold C$ has
dimension $1$.

When $a=1$, we assume $n\geq (r+1)/2$, and so we may apply  Larsen
Theorem [L]. With the same argument developed in [A], Proposition
8, pg. 69, we deduce that the Picard group of $X$ is generated by
the hyperplane section. In particular $X$ is subcanonical. At this
point, Theorem 0.1 in the hypotheses (D) with $a=1$ is a
consequence of Corollary 4.3.

When $a\geq 2$, using Theorem 4.1 instead of Theorem 2.1, one
proves Theorem 0.1 in the hypotheses (D) using a similar argument
as in the proof of Theorem 0.1 in the hypotheses (B) with $a\geq
2$. $\square$

\bigskip
\bigskip

\noindent {\bf Remark 4.4.} As a further consequence of Theorem
4.1, we have that {\it{the set of smooth, projective, numerically
subcanonical and not of general type subvarieties of dimension $n$
in ${\bold P^{2n}}$, is bounded}}. $\square$
\bigskip
\bigskip
\bigskip

\noindent {\bf 5. The proof of Theorem 0.1 under the hypotheses
(E).}
\bigskip

We keep the notation  introduced in Section 4, and assume the
hypotheses (E). Fix any real number $\lambda$ such that
$$
0<\lambda <{\frac{(a-1)^{n-1}(2a+n-2)}{2n!t^{\frac{n}{k}}}}. \tag
5.1
$$
Let $X\subset V$ be a smooth subvariety of dimension $n$ and
degree $d$, and assume that
$$
h^0(X,\Cal O_X(aK_X-bH_X))-c\left( \sum_{h=1}^{n-1}
p_g(X^{(h)})\right)\leq \lambda d^{{\frac{n}{k}}+1}. \tag 5.2
$$
As in $(2.12)$, using Kawamata-Viehweg Vanishing Theorem, one
proves that, when $a\geq 2$,
$$
h^0(X,\Cal O_X(aK_X-bH_X))\geq {\frac{d}{2n!}}(a-1)^{n-1}(2a+n-2)
n_1^n+dO(n_1^{n-1}),
$$
where  $n_1=e-q$, $e=(2g-2-(n-1)d)/d$, and $K_V=qH_V$ in
$H^2(V,{\bold C})$. Therefore we have
$$\multline
h^0(X,\Cal O_X(aK_X-bH_X))-c\left( \sum_{h=1}^{n-1}
p_g(X^{(h)})\right)\\
\geq {\frac{d}{2n!}}(a-1)^{n-1}(2a+n-2)
n_1^n+dO(n_1^{n-1})-c\left( \sum_{h=1}^{n-1} p_g(X^{(h)})\right).
\endmultline \tag 5.3
$$
Now using a Barth-Lefschetz type of argument (see diagram
$(2.4)$), our hypotheses on the codimension $k$ of $X$ and on the
cohomology of $V$ imply that, for any $1\leq i\leq k$,
$$
c_i(N_{X,V})=n_iH_X^i,
$$
in $H^{2i}(X,{\bold C})$, for a suitable rational number $n_i$. By
Lemma $(4.2)$ we deduce that, for any $1\leq i\leq k$,
$$
|n_i|\leq O(n_1^i).
$$
Using Proposition 3.2 and arguing as in the proof of Lemma 3.4, it
follows that
$$
p_g(X^{(h)})\leq dO(n_1^h)
$$
for any $1\leq h\leq n$. Hence we have
$$
\left( \sum_{h=1}^{n-1} p_g(X^{(h)})\right)\leq dO(n_1^{n-1}).
$$
Therefore, from $(5.3)$, we get
$$\multline
h^0(X,\Cal O_X(aK_X-bH_X))-c\left( \sum_{h=1}^{n-1}
p_g(X^{(h)})\right)\geq \\
{\frac{d}{2n!}}(a-1)^{n-1}(2a+n-2) n_1^n+dO(n_1^{n-1}).
\endmultline \tag 5.4
$$
Now fix any real number $\mu$ such that $\lambda<\mu<\frac
{(a-1)^{n-1}(2a+n-2)}{2n!t^{\frac{n}{k}}}$ (compare with $(5.1)$),
and put
$$
\lambda_1={\root{n}\of{\frac{n!\mu}{2^{n-1}(a-1)^{n-1}(2a+n-2)}}}.
\tag 5.5
$$
We have $0<\lambda_1<\frac{1}{2\root{k}\of{t}}$ (see $(4.1.1)$).
Hence, from Theorem 4.1, we may assume
$n_1>2\lambda_1d^{\frac{1}{k}}$. From $(5.4)$ we obtain
$$
h^0(X,\Cal O_X(aK_X-bH_X))\geq
\frac{(a-1)^{n-1}}{n!}2^{n-1}(2a+n-2)\lambda_1^nd^{{\frac{n}{k}}+1}+O(d^{{\frac{n+k-1}{k}}}).
\tag 5.6
$$
In view of $(5.5)$, comparing $(5.6)$ with $(5.2)$ we deduce that
$d$ is bounded.

\smallskip
One proves the sharpness of the estimate in a similar way as in
the hypotheses (A). This concludes the proof of Theorem 0.1 under
the hypotheses (E). $\square$
\bigskip
\bigskip

{\bf Aknowledgements}: We would like to thank F. Zak for some
useful discussion on the subject of this paper.
\bigskip
\bigskip

{

\end